\documentclass[11pt]{amsart}
\usepackage{amssymb}
\usepackage{amsxtra}
\theoremstyle{plain}
\newcommand{\cleqn}{\setcounter{equation}{0}}
\newcommand{\clth}{\setcounter{theorem}{0}}
\newcommand {\sectionnew}[1]{\section{#1}\cleqn\clth}
\newtheorem{theorem}{Theorem}[section]
\newtheorem{lemma}[theorem]{Lemma}
\newtheorem{definition-theorem}[theorem]{Definition-Theorem}
\newtheorem{proposition}[theorem]{Proposition}
\newtheorem{corollary}[theorem]{Corollary}
\newtheorem{definition}[theorem]{Definition}
\newtheorem{example}[theorem]{Example}
\newtheorem{remark}[theorem]{Remark}
\newtheorem{notation}[theorem]{Notation}
\newcommand \bth[1] { \begin{theorem}\label{t#1} }
\newcommand \ble[1] { \begin{lemma}\label{l#1} }

\newcommand \bpr[1] { \begin{proposition}\label{p#1} }
\newcommand \bco[1] { \begin{corollary}\label{c#1} }
\newcommand \bde[1] { \begin{definition}\label{d#1}\rm }
\newcommand \bex[1] { \begin{example}\label{e#1}\rm }
\newcommand \bre[1] { \begin{remark}\label{r#1}\rm }

\newcommand \bnota[1] { \begin{notation}\label{n#1}\rm }
\renewcommand {\eth} { \end{theorem} }
\newcommand {\ele} { \end{lemma} }

\newcommand {\epr} { \end{proposition} }
\newcommand {\eco} { \end{corollary} }
\newcommand {\ede} { \end{definition} }
\newcommand {\eex} { \end{example} }
\newcommand {\ere} { \end{remark} }

\newcommand {\enota} { \end{notation} }
\newcommand \thref[1]{Theorem \ref{t#1}}

\newcommand \prref[1]{Proposition \ref{p#1}}

\newcommand \lb[1]{\label{#1}}
\def \Rset {{\mathbb R}}         %mathsets
\def \Cset {{\mathbb C}}

\def \De {\Delta}   % Greek letters

\def \al {\alpha}

\def \Ga {\Gamma}

\def \ci  {\circ}

\def \Ad { {\mathrm{Ad}} }

\def \Lie { {\mathrm{Lie \,}} }
\def \gl  {\mathfrak{gl}}

\def \GL  {\mathrm{GL}}

\DeclareMathOperator \Gr { {\mathrm{Gr}} }

\DeclareMathOperator \Vect { {\mathrm{Vect}} }

\begin{document}
%%%%%%%%%%%%%%%%%%%%%%%%%%%%%%%%%%%%%%%%%%%%%%%%%%%%%%%%%%%%%%%%%%%%%%%%%%%
%%%%%%%%%%%%%%%%%%%%%%    Title    %%%%%%%%%%%%%%%%%%%%%%%%%%%%%%%%%%%%%%%%
\title[Cyclicity of Lusztig's stratification of Grassmannians]
{Cyclicity of Lusztig's stratification of Grassmannians and Poisson geometry}
\author[Milen Yakimov]{Milen Yakimov}
\address{ Department of Mathematics \\
Louisiana State Univerity \\
Baton Rouge, LA 70803 and
Department of Mathematics \\
University of California \\
Santa Barbara, CA 93106 \\
U.S.A.
}
\email{yakimov@math.lsu.edu}
\date{}
\begin{abstract} 
We prove that the standard Poisson structure
on the Grassmannian $\Gr(k, n)$ is invariant under the 
action of the Coxeter element $c =(1 2 \ldots n)$. In particular, 
its symplectic foliation is invariant under $c$. As a corollary,
we obtain a second, Poisson geometric proof of the result
of Knutson, Lam, and Speyer that the Coxeter element $c$
interchanges the Lusztig strata of $\Gr(k, n)$.
We also relate the main result to known anti-invariance properties 
of the standard Poisson structures on cominuscule flag varieties.  
\end{abstract}
\maketitle
%%%%%%%%%%%%%%%%%%%%   Introduction   %%%%%%%%%%%%%%%%%%%%%%%%%%%%%%%%%%%%%%%%
\sectionnew{Introduction}\lb{Introduction}
For the purpose of the study of canonical bases, Lusztig defined \cite{L} 
the totally nonnegative part $(G/P)_{\geq 0}$ of an arbitrary 
complex flag variety $G/P$. He also constructed an algebro-geometric 
stratification of $G/P$ and conjectured
that intersecting this stratification with $(G/P)_{\geq 0}$
is producing a cell decomposition of $(G/P)_{\geq 0}$.
This was latter proved by Rietsch 
in \cite{R}. Both the non-negative part $(G/P)_{\geq 0}$
and the Lusztig stratification of a flag variety 
were studied in recent years from many different combinatorial 
and Lie theoretic points of view. 

In a recent work Knutson, Lam, and Speyer proved
that the Lusztig stratification of the Grassmannian 
$\Gr(k,n)$ has a remarkable cyclicity property. If 
$c$ denotes the Coxeter element $(1 2 \ldots n)$ of $S_n$  
and the permutation matrix in $\GL_n(\Cset)$ 
which represents it, then $c$ permutes the 
strata of the Lusztig stratification of $\Gr(k,n)$.

In this note we give a Poisson geometric proof of this fact.
We also prove a stronger invariance property of a 
Poisson structure on $\Gr(k,n)$. In \cite{GY}, jointly with Goodearl, 
we found a Poisson geometric interpretation of 
the Lusztig stratification of any flag variety $G/P$.
For a choice of opposite Borel subgroups $B$ and $B^-$
of $G$ such that $B \subset P$ one defines the standard 
Poisson structure $\pi_{G/P}$ on $G/P$ which is 
invariant under the action of the maximal torus
$T = B \cap B^-$, see \cite{GY} for details.
According to \cite[Theorem 0.4]{GY}
the $T$-orbits of symplectic leaves of $\pi_{G/P}$ are
exactly the Lusztig strata.     

In the case of the complex Grassmannian $\Gr(k, n)$
the standard Poisson structure is given by
\[
\pi_{k,n} = - \sum_{1\leq i < j \leq n} 
\chi(E_{ij}) \wedge \chi(E_{ji})
\]
where $\chi \colon \gl_n(\Cset) \to \Vect(\Gr(k, n))$ denotes
the induced infinitesimal action from the left action
of $\GL_n(\Cset)$ on $\Gr(k, n)$ and $E_{ij}$ denote the 
elementary matrices. This Poisson structure 
is invariant under the action of the maximal torus $T_n$
of diagonal matrices in $\GL_n(\Cset)$. For each $w \in S_n$
denote by the same letter the corresponding permutation 
matrix in $\GL_n(\Cset)$. As before $c$ denotes 
the permutation matrix corresponding to the Coxeter element
$(1 2 \ldots n)$. The main result of this paper
is:

\bth{inv} Multiplication by $c$ is a Poisson 
automorphism of $(\Gr(k, n),$ $\pi_{k,n})$.
\eth

It is well known that the action of the permutation matrix $w_\ci$
corresponding to the longest element of $S_n$ is an anti-Poisson 
automorphism, see Section 2 for details. Thus \thref{inv} implies:

\bco{Dihedr} The actions of $w_\ci$, $c$, and $T_n$ generate 
an action of $I_2(n) \ltimes T_n$ by Poisson and anti-Poisson 
automorphisms of $(\Gr(k,n), \pi_{k,n})$ 
where $I_2(n)$ denotes the dihedral group of order $2n$.
\eco

The Lusztig stratification of the Grassmannian 
$\Gr(k, n)$ is defined as follows, see \cite{L} for details. 
Let $B$ and $B_-$ be 
the standard Borel subgroups of $\GL_n(\Cset)$ 
consisting of upper and lower triangular matrices.
Denote the maximal parabolic subgroup
\[
P_{k,n}= \{ \left( \begin{smallmatrix} a & b \\ 0 & c 
\end{smallmatrix} \right) \in \GL_{n}(\Cset)
\mid 
a \in M_{k, k}, b \in M_{k,n-k}, 
c \in M_{n-k,n-k} \}
\]
of $\GL_n(\Cset)$ and the induced map
\[
q \colon \GL_n(\Cset)/B \to \GL_n(\Cset)/P_{k,n} \cong
\Gr(k,n).
\]
The strata in the Lusztig stratification of 
$\Gr(k,n)$ are given by 
\[
R_{v,w}= q( B_- \cdot v B \cap B \cdot w B), \quad v 
\in (S_n)^{S_k \times S_{n-k} }_{\mathrm{max}}, 
w \in S_n, v \leq w.
\]
Here $\leq$ refers to the Bruhat order. We denote by
$S_k \times S_{n-k}$ the subgroup of $S_n$ consisting of
of those $u \in S_n$ such that $u(i) \leq k$ for 
$i \leq k$ and $u(i) \geq k+1$ for $i \geq k+1$.
Finally, $(S_n)^{S_k \times S_{n-k}}_{\mathrm{max}}$
denotes the set of maximal length representatives 
for the cosets $S_n/(S_k \times S_{n-k})$

The symplectic foliation of a Poisson structure is uniquely determined 
by it. Thus the $T_n$-orbits of leaves of $\pi_{k,n}$ 
(which are exactly the Lusztig strata) are an invariant 
of the pair ($\pi_{k,n}$, $T_n$-action). 
Therefore \thref{inv} gives 
a second proof of the result of Knutson, Lam, and Speyer 
that the action of the Coxeter element $c$ on $\Gr(k,n)$ interchanges 
the Lusztig strata. In fact \thref{inv} is equivalent 
to the stronger statement:

{\em{The action of $c$ on $\Gr(k,n)$ restricts to Poisson isomorphisms
between various Lusztig strata $(R_{v,w}, \pi_{k,n}|_{R_{v,w}})$
considered as regular Poisson varieties.}}

Finally we trace the roots of this phenomenon from a Poisson 
geometric point of view. It is well known that on any flag variety $G/P$  
the standard Poisson structure $\pi_{G/P}$ is anti-invariant under the action 
of any representative $\dot{w}_\ci$ of the longest element of 
the Weyl group $W$ of $G$. If, in addition, $P$ is cominuscule, 
then \cite[Proposition 4.2]{GY} implies that the standard Poisson
structure on $G/P$ is anti-invariant under the action 
of any representative $\dot{w}_\ci^P$ of the longest element 
of the corresponding parabolic subgroup of $W$. In the special 
case of the Grassmannian the specific Coxeter 
element $c$ happens to be a $k$-th root 
of $\dot{w}_\ci^P \dot{w}_\ci$. Thus \thref{inv} claiming that 
the standard Poisson structure on $\Gr(k,n)$ is invariant under 
$c$ is a strengthening of \cite[Proposition 4.2]{GY}.
See Section 2 for more details. We do not know of 
good Poisson properties of roots of $\dot{w}_\ci^P \dot{w}_\ci$ 
for any other cominuscule flag varieties.
\\ \hfill \\
{\bf Acknowledgements.} The author is grateful to 
Allen Knutson, Thomas Lam, and David Speyer for sharing 
the results of their preprint \cite{KLS} with him, which 
inspired this work. We would like thank Ken Goodearl
whose numerous comments helped us to improve the exposition.
We would also like to thank the organizers of the conference 
on Noncommutative Structures in Mathematics and Physics for 
the opportunity to participate at this very interesting meeting.
The author's research was partially supported
by NSF grant DMS-0701107.
%%%%%%%%%%%%%%%%%%%%%%%%%%%%%%%%%%%%%%%%%%
\noindent
\sectionnew{Proof of \thref{inv}}
\noindent
{\em{Proof of \thref{inv}.}} The statement is equivalent to 
showing that 
\[
\sum_{1\leq i < j \leq n} 
\chi(\Ad_c(E_{ij})) \wedge \chi(\Ad_c(E_{ji}))
- \sum_{1\leq i < j \leq n} \chi(E_{ij}) \wedge 
\chi(E_{ji})=0;
\]
that is
\begin{equation}
\label{eq}
V:= \sum_{i=2}^n \chi(E_{1i}) \wedge \chi(E_{i1}) =0.
\end{equation}
We will check this on the open Schubert cell 
$B_- \cdot P_{k,n} \subset Gr(k,n)$. Since $V$ is an algebraic bivector 
field, this will establish \eqref{eq}.
Identify
\begin{equation}
\label{ident}
M_{n-k, k} \cong B^- \cdot P_{k,n} \subset \Gr(k,n), 
\quad X \mapsto \left( \begin{smallmatrix} I_k & 0 \\ X & I_{n-k} 
\end{smallmatrix} 
\right) \cdot P_{k,n}
\end{equation}
where $M_{n-k, k}$ denotes the space of $(n-k) \times k $ complex
matrices. Applying \cite[eq. (3.17)]{BGY} 
we get that under \eqref{ident} 
\begin{equation}
\label{eq1}
\chi(E_{1,i+k}) \mapsto - \sum_{p=1}^{n-k}\sum_{q=1}^k x_{p1} x_{iq}
\frac{\partial}{\partial x_{pq}}, \quad \mbox{for}
\; i =1, \ldots, n-k.
\end{equation}
It is obvious that 
\begin{equation}
\label{eq2}
\chi(E_{i+k, 1}) \mapsto \frac{\partial}{\partial x_{i1}}, \quad \mbox{for}
\; i =1, \ldots, n-k.
\end{equation}
Let $1 \leq i, j \leq k$. Then 
\[
\Ad_{\exp(sE_{ij})} 
\left( \begin{smallmatrix} I_k & 0 \\ X & I_{n-k} 
\end{smallmatrix} \right) =
\left( \begin{smallmatrix} I_k & 0 \\ 
X - s \sum_{p=1}^{n-k} X_{pi}E_{pj} & I_{n-k} 
\end{smallmatrix} 
\right),
\]
which implies that under \eqref{ident}
\begin{equation}
\label{eq3}
\chi(E_{i,j}) \mapsto - \sum_{p=1}^{n-k} x_{pi} \frac{\partial}{\partial x_{pj}}.
\end{equation}
The summation in \eqref{eq} can be taken from 
$i=1$ since $\chi(E_{1,1})\wedge \chi(E_{1,1})=0$.
Applying \eqref{eq1}, \eqref{eq2}, and \eqref{eq3} we obtain that under the 
identification \eqref{ident}
\begin{multline}
V|_{B_- \cdot P_{k,n}} \mapsto
\sum_{i=1}^k \sum_{p=1}^{n-k} \sum_{q=1}^{n-k}
x_{p1} x_{qi} \frac{\partial}{\partial x_{pi}}
\wedge \frac{\partial}{\partial x_{q1}} \\
- \sum_{i=1}^{n-k} \sum_{p=1}^{n-k} \sum_{q=1}^{k}
x_{p1} x_{iq} \frac{\partial}{\partial x_{pq}}
\wedge \frac{\partial}{\partial x_{i1}} = 0.
\end{multline}
This implies \eqref{eq} and the statement of the 
Theorem.
\qed

For an arbitrary complex simple group $G$ and a maximal parabolic
subgroup $P$, one defines the standard Poisson structure 
\begin{equation}
\label{GP}
\pi_{G/P} = - \chi(r_G)
\end{equation}
on the flag variety $G/P$ induced from a compatible triangular 
decomposition of $G$ (a pair of Borel subgroups $B$ and $B_-$, such 
that $B \cap B_- =T$ is a maximal torus of $G$ and $B \subset P$), 
see e.g. \cite{GY}. Here $\chi \colon \wedge^2 \Lie(G) \to 
\Ga(T G/P, G/P)$ denotes the induced action from the infinitesimal 
action of $\Lie(G)$. The standard $r$-matrix $r_G \in \wedge^2 \Lie(G)$ 
obtained from the triangular decomposition of $G$ is given by:
\[
r_G = \sum_{\al \in \De_+} e_\al \wedge f_\al
\]
where $e_\al$ and $f_\al$ are appropriately normalized
root vectors of $\Lie(G)$ and $\De_+$ is the 
set of positive roots of $G$, cf. \cite{GY}.  
It is obvious that the action of any representative 
$\dot{w}_\ci$ of the longest element of the Weyl group $W$ of $G$ 
on $(G/P, \pi_{G/P})$ is anti-Poisson, since $\Ad_{\dot{w}_\ci}$
interchanges $e_\al$ and $f_\al$.

Denote the Levi factor of $P$ containing $T$ by $L$, and the 
longest element of the subgroup of $W$ corresponding to $L$ by 
$w_\ci^P$. Let $\dot{w}_\ci^P$ be any representative of 
$w_\ci^P$ in the normalizer of $T$.

Recall that among several equivalent definitions/characterizations
of cominuscule parabolic subgroups: a parabolic subgroup $P$ of $G$ is
cominuscule if and only if its unipotent radical is abelian.
According to \cite[Proposition 4.2]{GY}, if $P$ cominuscule,
then $\pi_{G/P}$ is also given by
\[
\pi_{G/P} = - \chi(r_L).
\] 
where $r_L \in \wedge^2 \Lie(L)$ is the standard $r$-matrix of $L$.
Thus the action of $\dot{w}_\ci^P$ on 
$(G/P, \pi_{G/P})$ is anti-Poisson as well. So:

\bpr{1}
For any cominuscule parabolic subgroup $P$
of a complex simple Lie group $G$, the action of 
$\dot{w}_\ci^P \dot{w}_\ci$ on $(G/P, \pi_{G/P})$ 
is Poisson.
\epr

In the special case of the Grassmannian $\Gr(k,n)$ 
\[
w_\ci^P w_\ci = c^k
\]
for the particular Coxeter element $c$.
Taking powers of this product, we see that the action of 
$c^{gcd(k,n)}$ on $\Gr(k,n)$ is Poisson. In the case 
when $k$ and $n$ are relatively prime this gives yet 
another proof of \thref{inv}. One could argue that  
\thref{inv} holds because it is true for relatively prime 
$k$ and $n$, and its statement (cf. also its proof) is independent
of the numerical properties of $k$ and $n$. 

Conceptually the invariance of $\pi_{k,n}$ under 
$c$ is the result of a two step process: 

1. From \prref{1} one has the invariance of $\pi_{G/P}$ under 
the product $\dot{w}_\ci^P \dot{w}_\ci$ of the longest elements 
of the Weyl groups of $G$ 
and the Levi subgroup $L$, for arbitrary cominuscule
flag variety $G/P$.

2. In the case of the Grassmannian the Coxeter element 
$c$ which is a $k$-th root of $\dot{w}_\ci^P \dot{w}_\ci$ acts 
by Poisson automorphisms of $(\Gr(k,n), \pi_{k,n})$ as well.

The special property of the Coxeter element $c=(1 2 \ldots n)$ is that 
the other Coxeter elements of $S_n$ are not roots of 
$\dot{w}_\ci^P \dot{w}_\ci$. 
%%%%%%%%%%%%%%%%%%%%%%%%%%%%%%%%%%%%%%%%%%
\sectionnew{Corollaries}
The symplectic foliation of a Poisson manifold $(M, \pi)$ is an invariant 
of it. Similarly if a group $H$ act on $(M, \pi)$ by Poisson automorphisms
the partition of $M$ into $H$-orbits of symplectic leaves is an 
invariant of $(M, \pi)$ considered as a Poisson $H$-space. 
Since the partition of $\Gr(k,n)$ into
$T_n$-orbits of symplectic leaves of $(\Gr(k,n), \pi_{k,n})$ is exactly the 
Lusztig stratification of $\Gr(k,n)$ due to \cite[Theorem 0.4]{GY}
and $c$ normalizes $T_n$, \thref{inv} implies the following Theorem 
of Knutson, Lam, and Speyer:

\bth{Lus} The action of the permutation matrix corresponding to 
the Coxeter element $c= (1 2 \ldots n)$ on $\Gr(k,n)$
permutes the strata 
$R_{v,w}= q( B_- \cdot v B \cap B \cdot w B)$ 
of the Lusztig stratification.
\eth

As we pointed out, in addition, when $c$ maps one Lusztig stratum $R_{v_1,w_1}$
to another $R_{v_2,w_2}$, it matches the regular Poisson structures 
$\pi_{k,n}|_{R_{v_1,w_1}}$ and $\pi_{k,n}|_{R_{v_2,w_2}}$.

Finally let us point out that all constructions and invariance properties are valid 
over the reals since all constructions are derived from the real split form 
$\GL_n(\Rset)$.
%%%%%%%%%%%%%%%%%%%%%% References %%%%%%%%%%%%%%%%%%%%%%%%%%%%%%%%%%%%%%%

%%%%%%%%%%%%%%%%%%%%%%%%%%%%%%%%%%%%%%%%%%%%%%%%%%%%%%%%%%%%%%%%%%%%%%%%%%%%%%%
%%%%%%%%%%%%%%%%%%%%%%%%%%%%%%%%%%%%%%%%%%%%%%%%%%%%%%%%%%%%%%%%%%%%%%%%%%%%%%
\end{document}